\documentclass[12pt]{article}
\newenvironment{algorithm}[1]
{
\begin{figure}[h]
  \begin{center}
    {\bf Algoritm: #1}\\*
    \begin{tabular}{|p{353pt}|} \hline  %
}
{
 \\ \hline
 \end{tabular}
 \end{center}
\end{figure}
}
\def\Z{{\mathchoice {\hbox{$\textstyle\sf Z\kern-0.4em Z$}}
{\hbox{$\textstyle\sf Z\kern-0.4em Z$}}
{\hbox{$\scriptstyle\sf Z\kern-0.3em Z$}}
{\hbox{$\scriptscriptstyle\sf Z\kern-0.2em Z$}}}}
\def\Q{{\mathchoice {\setbox0=\hbox{$\displaystyle\rm
Q$}\hbox{\raise
0.15\ht0\hbox to0pt{\kern0.4\wd0\vrule height0.8\ht0\hss}\box0}}
{\setbox0=\hbox{$\textstyle\rm Q$}\hbox{\raise
0.15\ht0\hbox to0pt{\kern0.4\wd0\vrule height0.8\ht0\hss}\box0}}
{\setbox0=\hbox{$\scriptstyle\rm Q$}\hbox{\raise
0.15\ht0\hbox to0pt{\kern0.4\wd0\vrule height0.7\ht0\hss}\box0}}
{\setbox0=\hbox{$\scriptscriptstyle\rm Q$}\hbox{\raise
0.15\ht0\hbox to0pt{\kern0.4\wd0\vrule height0.7\ht0\hss}\box0}}}}

\newcommand{\binom}[2]{\Bigl(\!\!
\begin{tabular}{c}
  $#1$ \\
  $#2$
\end{tabular}\!\!\Bigr)}
\newcommand{\Ti}[1]{$\!\!$#1$\!\!$&}
\newcommand{\Tie}[1]{$\!\!$#1$\!\!$\\ \hline}
\newcommand{\Tm}[1]{$\!\!{#1}\!\!$&}
\newcommand{\Tvoid}{&}
\newcommand{\Td}[3]{$\!\!\!$\begin{tabular}{c}
       \\[-11pt]
  {#1} \\
  {#2} \\
  {#3} \\
\end{tabular}$\!\!\!$&}
\newcommand{\Tde}[3]{$\!\!$\begin{tabular}{c}
       \\[-11pt]
  {#1} \\
  {#2} \\
  {#3} \\
\end{tabular}$\!\!$\\ \hline}
\newcommand{\Tc}[4]{$\!\!\!\!$\begin{tabular}{c}
       \\[-11pt]
  {#1} \\
  {~~#2} $\!\bullet$\\ 
  {#3} \\
\end{tabular}$\!\!\!$&}
\newcommand{\Tce}[4]{$\!\!\!$\begin{tabular}{c}
       \\[-11pt]
  {#1} \\
  {~~#2}  $\!\bullet$\\ 
  {#3} \\
\end{tabular}$\!\!\!$\\ \hline}

\begin{document}
\begin{center}
\begin{center}
{\Large\bf Computation of Cohomology of Lie Algebra of Hamiltonian Vector Fields
by Splitting Cochain Complex into Minimal Subcomplexes}
\end{center}
\begin{center}
{\large  V. V. Kornyak}
\end{center}
\begin{center}
\emph{Laboratory of Information Technologies\\
Joint Institute for Nuclear Research\\
141980 Dubna, Russia \\
fax: (+7)-(09621)-65145\\
email: kornyak@jinr.ru}
\end{center}
\end{center}
{\small
Computation of homology or cohomology is intrinsically a problem of high combinatorial
complexity.
Recently we proposed a new efficient algorithm for
computing cohomologies of Lie algebras and superalgebras. This
algorithm is based on partition of the full cochain complex into
minimal subcomplexes. The algorithm was implemented as a C
program {\bf LieCohomology}. In this paper we present results of applying the program
{\bf LieCohomology} to the algebra of hamiltonian vector fields $\mathrm{H}(2|0)$. We
demonstrate that the new approach is much more efficient comparing
with the straightforward one. In particular, our computation reveals some new
cohomological classes for the algebra $\mathrm{H}(2|0)$ (and also for the Poisson algebra $\mathrm{Po}(2|0)$).
}

\section{Introduction}
\noindent
Cohomology is defined by \emph{cochain complex}
\begin{equation}
0\to C^0\stackrel{d^0}{\longrightarrow}\cdots
\stackrel{d^{k-2}}{\longrightarrow}C^{k-1}\stackrel{d^{k-1}}
{\longrightarrow} C^k\stackrel{d^k}{\longrightarrow}C^{k+1}
\stackrel{d^{k+1}}{\longrightarrow}\cdots.
\label{cocomplex}
\end{equation}
Here $C^k$ are linear spaces (more generally, abelian groups), graded by the
integer number $k,$ called \emph{dimension} or \emph{degree} (depending on the context).
The elements of the spaces $C^k$ are called \emph{cochains}.

The linear mappings $d^k$ are called \emph{differentials} (or \emph{coboundary operators}).
The main property of these mappings is ``their squares are equal to
zero'': $d^k \circ d^{k-1} = 0$.

The elements of the space $Z^k = \mathrm{Ker}\ d^k$ are called \emph{cocycles}.
The elements of the space $B^k = \mathrm{Im}\ d^{k-1}$ are called \emph{coboundaries}.
Note that $B^k \subseteq Z^k.$

The $k$th \emph{cohomology} is the quotient space
$$H^k = Z^k/B^k \equiv {\mathrm{Ker}\ d^k}/{\mathrm{Im}\ d^{k-1}}.$$

There are many cohomological theories designed for
investigation of different mathematical structures and the space $H^k$ carries important information
about peculiarities in these structures. The only difference between
cohomological theories lies in the constructions of the cochain spaces and coboundary
operator. These constructions depend on the underlying mathematical
structures.

The cohomology of the Lie (super)algebra $A$ in the module $X$ is defined
via cochain complex (\ref{cocomplex}) in which (see, e.g., \cite{Fuks})
the cochain spaces
$C^k = C^k(A;X)$ consist of  super skew-symmetric
$k$-linear
mappings $A \times \cdots \times A \to X,$  $C^0 = X$ by definition.
Super skew-symmetry means symmetry with respect to swapping of two adjacent
\emph{odd} cochain arguments
and antisymmetry  for any other combination of parities for adjacent pair.

The differential $d^k$  takes the form\footnote{This version of formula for differential corresponds to the algorithm
used in the program {\bf LieCohomology}.}
\begin{eqnarray}
(d^kc)(a_0,\ldots,a_k)& = &-\sum_{0\leq i<j\leq k}(-1)^{s(a_i)+s(a_j)+p(a_i)p(a_j)}
c([a_i,a_j],a_o,\ldots,\widehat{a_i},\ldots,\widehat{a_j},\ldots,a_k)\nonumber\\
& & -\sum_{0\leq i \leq
k}(-1)^{s(a_i)}a_ic(a_o,\ldots,\widehat{a_i},\ldots,a_k),
\label{diff}
\end{eqnarray}
where the functions $c(\ldots)$ are elements of cochain spaces; $a_i \in A$; $p(a_i)$
is the parity of  $a_i$; $s(a_i) = i,$ if $a_i$ is even element and $s(a_i)$
is equal to the number of even elements in the sequence $a_0,\ldots, a_{i-1}$,
if $a_i$ is odd element.
In the case of trivial module (i.e., if $ax = 0$ for all
$a \in A$ and  $x \in X$) one uses as a rule the notation $H^k(A).$

In papers~\cite{KornProg99,Kornyak99b,Kornyak00Co1,KornProg01} we presented an algorithm
for computation of Lie (super)algebra cohomologies.
These papers contain also description of its C
implementation and some results obtained with the help of codes designed.
This algorithm computes cohomology of Lie (super)algebra $A$ over module $X$
in a straightforward way, i.e.,
for cochain complex (\ref{cocomplex})
the algorithm constructs the full set of basis super skew-symmetric monomials forming the space $C^{k}$,
generates subsequently all basis monomials in the space $C^{k+1}$,
computes the differentials corresponding to these monomials to obtain the set of
linear equations determining the space
of cocycles
\begin{equation}
Z^k = \mathrm{Ker}\ d^k = \{C^k \; | \; dC^k = 0\},
\end{equation}
constructs the space of coboundaries
\begin{equation}
B^k = \mathrm{Im}\ d^{k-1} = \{C^k \; | \; C^k = dC^{k-1}\}.
\end{equation}
Finally, the algorithm constructs the basis elements of quotient space
\begin{equation}
H^k(A;X) = Z^k/B^k.
\end{equation}
This last step is based on the Gauss elimination procedure.

The main difficulty in computing cohomology results from the very high dimensions
of the spaces $C^k$: for $n$-dimensional ordinary Lie algebra and $p$-dimensional
module $$\dim C^k = p\, \binom{n}{k},$$ and for $(n|m)$-dimensional Lie superalgebra
$$\dim C^k = p \sum_{i=0}^{k}\binom{n}{k-i}\binom{m+i-1}{i} \equiv p\, \binom{n}{k} + p \sum_{i=1}^{k}\binom{n}{k-i}\binom{m+i-1}{i}.$$

In many cases it is possibly to extract some easier to handle subcomplexes of the
full cochain complex (\ref{cocomplex}). The partition of cochain complex for a graded
algebra and module into homogeneous components is a typical example.
In many papers (see, e. g., \cite{GKF,Perchik,GSH}) more special subcomplexes were
used successfully to obtain new results in the theory of cohomology of Lie
(super)algebras.\footnote{The main trick consists in imposing some restrictions
on the elements of $C^k$ and proving the invariance of these restrictions with respect
to the differential.}

The main idea of the new algorithm presented in \cite{KornCASC01,KornProg02,KornProg02a}
is to extract
the minimal possible subcomplexes from complex (\ref{cocomplex}) and to carry
computations within these subcomplexes. There are two versions of the algorithm.
One of them is applied when the cochain spaces under consideration are infinite-dimensional
(or their dimensions are too large to fit the available memory), but the minimal subcomplexes
contain finite-dimensional spaces of $k$-cochains. Another version of the algorithm is applied when
it is possible to construct the full space $C^k.$  Below we present this version in the pseudocode form.
\begin{algorithm}{ComputeCohomology\label{ComputeCohomology}}
\item[\bf{~Input:}]  $~A$, Lie~(super) algebra; $X$, module;\\
\hspace*{42pt} $k$, cohomology degree; $g$, grade
\item[\bf{~Output:}] $BH^k_g$, set of basis cohomological classes
\item[\bf{~Local:}] $~~~M^k_g$, full set of $k$-cochain monomials (basis of $C^k_g$);\\
\hspace*{42pt} $s$, current subcomplex: $C^{k-1}_{g,s}\stackrel{d^{k-1}_{g,s}}
{\longrightarrow} C^k_{g,s}\stackrel{d^k_{g,s}}{\longrightarrow}C^{k+1}_{g,s};$\\
\hspace*{42pt} $m^k_g \in M^k_g$, starting monomial for constructing subcomplex $s$;\\
\hspace*{42pt} $M^k_{g,s}$, set of $k$-cochain monomials involved in subcomplex $s$;\\
\hspace*{42pt} $BH^k_{g,s}$, set of basis cohomological classes in subcomplex $s$

~1: $BH^k_g:=\emptyset$\\
~2: $M^k_g:= {\bf GenerateMonomials}(A,\,X,\,k,\,g)$\\
~3: \bf{while} $M^k_g\neq \emptyset$ \bf{do}\\
~4: $~~~$ $m^k_g:={\bf ChooseMonomial}(M^k_g)$\\
~5: $~~~$ $\{s,\,M^k_{g,s}\}:={\bf ConstructSubcomplex}(m^k_g)$\\
~6: $~~~$ $BH^k_{g,s}:={\bf ComputeCohomologyInSubcomplex}(s)$\\
~7: $~~~$ \bf{if} $BH^k_{g,s}\neq\emptyset$ \bf{then}\\
~8: $~~~~~~$ $BH^k_g:=BH^k_g \cup BH^k_{g,s}$\\
~9: $~~~$ \bf{fi}\\
10: $~~~$ $M^k_g:=M^k_g\setminus M^k_{g,s}$\\
11: \bf{od}\\
12: \bf{return} $BH^k_g$\\
\end{algorithm}

\noindent
Here the subalgorithm {\bf GenerateMonomials} generates the full set $M^k_g$ of
super skew-symmetric monomials
$$
c(\alpha_{i_1},\ldots,\alpha_{i_k};\xi_\iota)\equiv c(\alpha_{i_1})\wedge\cdots\wedge c(\alpha_{i_k})
\otimes \xi_\iota\equiv \alpha'_{i_1}\wedge\cdots\wedge \alpha'_{i_k}\otimes \xi_\iota
$$
forming basis of the cochain space $C^k$ in the grade $g;$ $\alpha_{i_j} \in A$ and $\xi_\iota \in X$ are basis
elements of algebra and module; $\alpha'_i$ is the dual to $\alpha_i$ element. The subalgorithm {\bf ChooseMonomial}
takes some monomial $m^k_g \in M^k_g.$
This monomial is starting monomial for constructing the subcomplex $s$ by the subalgorithm {\bf ConstructSubcomplex}.
The subalgorithm {\bf ComputeCohomologyInSubcomplex} computes basis cohomological classes $BH^k_{g,s}$ in the subcomplex
$s$ by the straightforward algorithm described above.

\section{Computation of $H^k_g(\mathrm{H(2|0)})$}

In this section we present the results of computation of cohomology in the trivial module for Lie algebra
  $\mathrm{H(2|0)}$
of formal hamiltonian vector fields on the $2$-dimensional simplectic manifold. We
describe also cohomological classes up to grade 8 for the Poisson algebra $\mathrm{Po(2|0)}$ which
is a central extension of the algebra  $\mathrm{H(2|0)}$.

The hamiltonian algebra $\mathrm{H(2n|m)}$ is an algebra of vector fields (see, e.g., \cite{Leites}) acting on the $(2n|m)$ supermanifold
and preserving the following 2-form
$$\sum_{i=1}^n dp_i \wedge dq_i + \sum_{j=1}^m d\theta_j \wedge d\theta_j,$$
where $p_1,\ldots,p_n; q_1,\ldots,q_n$ and $\theta_1,\ldots\theta_m$ are even and odd local variables
on the supermanifold, respectively. The elements of $\mathrm{H(2n|m)}$ can be expressed in terms of
\emph{generating function} $f(p_1,\ldots,p_n; q_1,\ldots,q_n; \theta_1,\ldots\theta_m)$ by the formula
\begin{equation}
\sum_{i=1}^n \left(\frac{\partial f}{\partial p_i}\frac{\partial}{\partial q_i} -
                     \frac{\partial f}{\partial q_i}\frac{\partial}{\partial p_i}\right)
   - (-1)^{p(f)}\sum_{j=1}^m  \frac{\partial f}{\partial \theta_j}\frac{\partial}{\partial \theta_j},
   \label{vf}
\end{equation}
where $p(f)$ is parity of the function $f$ (this function is called usually \emph{hamiltonian}).
Thus one can consider the formal hamiltonian vector fields as linear combinations of monomials in the variables  $p_i,$ $q_i$ and
$\theta_j$ (except of the monomial 1). Considering these monomials as basis elements of $\mathrm{H(2n|m)}$
and using prescribed $\Z$-grading for the variables $p_i,$ $q_i$ and $\theta_j$ one can impose $\Z$-grading $\mathrm{gr}()$ on the algebra
$\mathrm{H(2n|m)}.$ The standard grading is $\mathrm{gr}(p_i) = \mathrm{gr}(q_i) = \mathrm{gr}(\theta_j) = 1.$
For the standard grading the grade of algebra element corresponding to some monomial is equal to
the grade of this monomial minus 2 (due to two differentiations in the terms of expression (\ref{vf})).

Since the hamiltonian algebra is very important in both classical and quantum physics many efforts were
applied to investigation of its cohomology. Most advanced results were obtained for the finite-dimensional algebras of the form
$\mathrm{H(0|m)}$ \cite{FuchsLeites84,Gruson97}. Nevertheless the hamiltonian algebras on supermanifolds with
nonzero even dimension are more important in applications but computation of their cohomologies is much more difficult task.
Some results about cohomologies of such algebras were obtained in \cite{GKF,Perchik,GSH}.
In the paper \cite{GKF} some elements of  $H^k_g(\mathrm{H(2|0)})$ were calculated by considering special subcomplexes
(and using computer partially). We present here all cohomological
classes (without discussing their meaning and interpretation) from $H^k_g(\mathrm{H(2|0)})$ up to grade 8.

The results of computation are summarized in Table \ref{resulttable}. The boxes of this table corresponding to cochain degree
$k$ and cochain grade $g$ contain the following information: $\dim C^k_g,$ dimension of the full space of $k$-cochains in grade $g$;
number of minimal subcomplexes
$C^{k-1}_{g,s}\stackrel{d^{k-1}_{g,s}} {\longrightarrow} C^k_{g,s}\stackrel{d^k_{g,s}}{\longrightarrow}C^{k+1}_{g,s}$
extracted by the algorithm from the full complex; $\max \dim  C^k_{g,s},$ maximum dimension of the subspace of $(k,g)$-cochains
among all subcomplexes. The empty box means that $\dim C^k_g = 0,$ i.e. the space of $(k,g)$-cochains is empty.
The boxes marked by the bullet $\bullet$ contain nontrivial 1-dimensional cohomological classes. For example,
the box corresponding to the pair $(k,g) = (7,8)$ tells that $\dim C^7_8 = 25488,$ number of subcomplexes is 21,
$\max \dim  C^7_{8,s} = 3148$ and $\dim H^7_8 = 1.$

More detailed information about computation in $(k,g) = (7,8)$ is given in Table \ref{structtable}.
\begin{table}[h!]
\caption{Subcomplex structure for $(k,g)=(7,8)$}
\begin{center}\label{structtable}
\begin{tabular}{|r|r|r|r|r|c|c|}
\hline
$\dim C^{k-1}_{g,s}\!\!$&$\dim C^k_{g,s}\!\!$&$\dim C^{k+1}_{g,s}\!\!$&$\dim Z^k_{g,s}\!\!$&$\dim B^k_{g,s}\!\!$&$\dim H^k_{g,s}$&repeated\\ \hline
           0        &        1       &         1          &         0      &        0       &    0           & 2 \\ \hline
          12        &       17       &        11          &         9      &        9       &    0           & 2 \\ \hline
          72        &       80       &        54          &        43      &       43       &    0           & 2 \\ \hline
         223        &      243       &       167          &       130      &      130       &    0           & 2 \\ \hline
         507        &      540       &       375          &       292      &      292       &    0           & 2 \\ \hline
         909        &      976       &       702          &       520      &      520       &    0           & 2 \\ \hline
        1406        &     1536       &      1120          &       813      &      813       &    0           & 2 \\ \hline
        1928        &     2117       &      1578          &      1114      &     1114       &    0           & 2 \\ \hline
        2382        &     2652       &      1992          &      1387      &     1387       &    0           & 2 \\ \hline
        2695        &     3008       &      2286          &      1568      &     1568       &    0           & 2 \\ \hline
        2806        &     3148       &      2391          &      1640      &     1639       &    1           & 1 \\ \hline
\end{tabular}
\end{center}
\end{table}
In this table the columns $\dim Z^k_{g,s},$ $\dim B^k_{g,s}$ and $\dim H^k_{g,s}$ contain dimensions of cocycle, coboundary and
cohomology spaces in subcomplexes, respectively. On can see that there are 10 pairs of subcomplexes with repeated structure and the only
single subcomplex containing nontrivial cohomological class.

The full set of nontrivial (1-dimensional) cohomological classes is: $H^2_{-2}, H^5_{-2}, H^7_0$ (computed earlier)
and $H^7_8, H^{10}_6$ (computed by the new program). As to the Poisson algebra $\mathrm{Po(2|0)},$ the part of its cohomological classes
up to grade 8 coincides with those for $\mathrm{H(2|0)}$ except of $H^2_{-2}$.\footnote{The cocycle $H^2_{-2}$ describes the central extension
of the algebra $\mathrm{H(2|0)}$ to  $\mathrm{Po(2|0)}.$} $H^k_{g\leq8}(\mathrm{Po(2|0)})$ contains also four additional classes: $H^6_{-4}, H^8_{-2}, H^8_6, H^{11}_6.$
But all these classes are multiplicative consequences of the classes $H^5_{-2}, H^7_0, H^7_8, H^{10}_6.$
These classes can be expressed in the form $H^k_g = H^{k-1}_{g+2}\wedge c({\cal Z})$ due to the general property \cite{Kornyak00Co1}
of cohomology of algebras containing central element ${\cal Z}.$

\begin{table}[h!]
\caption{Computation of $H^k_g$ for  $(k,g) \in [1,\ldots,\infty)\otimes[-2,\ldots,8]$}
\begin{center}\label{resulttable}
\begin{tabular}{c|c|c|c|c|c|c|c|c|c|c|c|}
\Tm{k\backslash g}
       \Ti{-2}        \Ti{-1}      \Ti{0}         \Ti{1}          \Ti{2}          \Ti{3}           \Ti{4}            \Ti{5}            \Ti{6}            \Ti{7}             \Tie{8}
\Ti{~1}\Tvoid         \Td{2}{2}{1} \Td{3}{3}{1}   \Td{4}{4}{1}    \Td{5}{5}{1}    \Td{6}{6}{1}     \Td{7}{7}{1}      \Td{8}{8}{1}      \Td{9}{9}{1}      \Td{10}{10}{1}     \Tde{11}{11}{1}
\Ti{~2}\Tc{1}{1}{1}{1}\Td{6}{4}{2} \Td{11}{5}{3}  \Td{22}{6}{5}   \Td{33}{7}{7}   \Td{52}{8}{9}    \Td{71}{9}{11}    \Td{100}{10}{14}  \Td{129}{11}{17}  \Td{170}{12}{20}   \Tde{211}{13}{23}
\Ti{~3}\Td{3}{3}{1}   \Td{10}{4}{3}\Td{30}{7}{8}  \Td{60}{8}{13}  \Td{116}{9}{22} \Td{200}{10}{34} \Td{326}{11}{52}  \Td{502}{12}{72}  \Td{752}{13}{100} \Td{1078}{14}{133} \Tde{1515}{15}{177}
\Ti{~4}\Td{3}{3}{1}   \Td{14}{6}{4}\Td{45}{7}{11} \Td{108}{8}{22} \Td{242}{11}{44}\Td{466}{12}{78} \Td{857}{13}{135} \Td{1468}{14}{210}\Td{2426}{15}{326}\Td{3820}{16}{478} \Tde{5872}{17}{698}
\Ti{~5}\Tc{1}{1}{1}{1}\Td{12}{6}{3}\Td{41}{7}{9}  \Td{128}{10}{25}\Td{315}{11}{59}\Td{706}{12}{117}\Td{1432}{13}{222}\Td{2748}{16}{391}\Td{4949}{17}{671}\Td{8568}{18}{1078}\Tde{14240}{19}{1710}
\Ti{~6}\Tvoid         \Td{4}{4}{1} \Td{23}{7}{5}  \Td{90}{10}{18} \Td{264}{11}{50}\Td{688}{12}{114}\Td{1580}{15}{246}\Td{3382}{16}{483}\Td{6734}{17}{916}\Td{12766}{18}{1619}\Tde{23074}{19}{2806}
\Ti{~7}\Tvoid         \Tvoid       \Tc{6}{5}{2}{1}\Td{32}{8}{7}   \Td{135}{11}{25}\Td{412}{12}{70} \Td{1128}{15}{180}\Td{2730}{16}{396}\Td{6132}{17}{842}\Td{12818}{18}{1649}\Tce{25488}{21}{3148}{1}
\Ti{~8}\Tvoid         \Tvoid       \Tvoid         \Td{4}{4}{1}    \Td{33}{9}{7}   \Td{138}{10}{25} \Td{479}{13}{79}  \Td{1388}{16}{207}\Td{3606}{17}{510}\Td{8546}{18}{1125}\Tde{18963}{21}{2391}
\Ti{~9}\Tvoid         \Tvoid       \Tvoid         \Tvoid          \Td{1}{1}{1}    \Td{20}{8}{4}    \Td{99}{11}{17}   \Td{396}{14}{62}  \Td{1260}{15}{188}\Td{3576}{18}{489} \Tde{9104}{19}{1188}
\Ti{10}\Tvoid         \Tvoid       \Tvoid         \Tvoid          \Tvoid          \Tvoid           \Td{5}{5}{1}      \Td{46}{10}{8}    \Td{217}{13}{35}  \Td{818}{16}{120}  \Tce{2578}{17}{358}{1}
\Ti{11}\Tvoid         \Tvoid       \Tvoid         \Tvoid          \Tvoid          \Tvoid           \Tvoid            \Tvoid            \Td{10}{7}{2}     \Td{70}{10}{12}    \Tde{350}{15}{54}
\Ti{12}\Tvoid         \Tvoid       \Tvoid         \Tvoid          \Tvoid          \Tvoid           \Tvoid            \Tvoid            \Tvoid            \Tvoid             \Tde{10}{7}{2}
\end{tabular}
\end{center}
\end{table}

\section{Conclusion}
Our new algorithm demonstrates substantially higher efficiency in comparison with the old one.
For example, the program described in \cite{KornProg99} computes the case $(k,g) = (6,5)$ in 35 min 45 sec = 2145 sec
whereas the new program takes 54 sec for this task. For both runs we used PC Pentium III, 667MHz, 256MB RAM.
The superiority of the new program grows
with increasing of the task complexity.
Nevertheless, due to rapidly increasing computational complexity the presented results are not sufficient
to derive any general idea about the structure of cohomology ring  $H^*_*(\mathrm{H(2|0)})$.
Our computation was carried over the field of rational numbers $\Q.$
As profiling shows, the most time consuming part of computation by the program {\bf LieCohomology} is  multiprecision arithmetic.
This is common difficulty for almost all problems in computer algebra.
It seems that carrying computation over the finite fields, say $\Z_p,$ we can go to the grade 40-50 for the problem considered here,
but the results obtained in this way can be considered merely as hints.

\section*{Acknowledgements}
This work was partially supported by the grants RFBR 01-01-00708,
RFBF 00-15-96691 and INTAS 99-1222.


\end{document}